\documentclass[11pt]{article}
\usepackage{amsmath}
\usepackage{amsfonts}
\usepackage{amssymb}
\usepackage{euscript}
\newenvironment{proof}{\noindent {\it Proof.~~}\ }{\  \rule{1mm}{2mm}\medskip}

\newtheorem{theorem}{Theorem}
\newtheorem{lemma}[theorem]{Lemma}
\newtheorem{corollary}[theorem]{Corollary}
\newtheorem{proposition}[theorem]{Proposition}

\def\Z{\mathbb Z}

\def\G{\partial}

\newcommand{\subgp}[1]{\langle{#1}\rangle}

\newcommand{\ot}{\overline{T}}

\newcommand{\oa}{\overline{A}}

\newcommand{\oaa}{\overline{A^{-1}}}
\begin{document}
\title{On Group  bijections $\phi $  with $\phi(B)=A$ and $\forall a\in B, a\phi(a)  \notin A$ }
\author{ {Y. O. Hamidoune}\thanks{
UPMC Univ Paris 06,
 E. Combinatoire, Case 189, 4 Place Jussieu,
75005 Paris, France.}
}


\maketitle

\begin{abstract}

A {\em Wakeford pairing} from $S$ onto $T$  is a bijection $\phi : S \rightarrow T$ such that
$x\phi(x)\notin T,$ for every $x\in S.$ The number of such pairings will be denoted by $\mu(S,T)$.

 Let $A$ and $ B$ be  finite  subsets of  a   group $G$  with $1\notin B$ and $|A|=|B|.$  Also assume
 that the order of every element of $B$ is $\ge |B|$.

 Extending results due to Losonczy and Eliahou-Lecouvey,
 we show that $\mu(B,A)\neq 0.$

 Moreover we show that  $\mu(B,A)\ge  \min \{\frac{||B|+1}{3},\frac{|B|(q-|B|-1)}{2q-|B|-4}\},$
 unless there is $a\in A$ such that $|Aa^{-1}\cap B|=|B|-1$ or $Aa^{-1}$ is a progression.

 In particular,
either $\mu(B,B)
 \ge \min \{\frac{||B|+1}{3},\frac{|B|(q-|B|-1)}{2q-|B|-4}\},$ or
  for some $a\in B,$ $Ba^{-1}$ is a progression.

\end{abstract}

\hspace{1cm}{\bf MSC Classification:} 11B60, 11B34, 20D60.

\section{Introduction}

Let $A$ and $B$ be finite subsets of a group $G$ with and $|A|=|B|.$
A {\em Wakeford pairing} from $B$ onto $A$  is a bijection $\phi : B \rightarrow A$ such that
$x\phi(x)\notin A,$ for every $x\in B.$

Our pairings are dual to pairings used in litterature \cite{fl,los,el}. The two notions are equivalent up to replacing the group
its opposite group or by replacing $(A,B)$ by $(A^{-1},B^{-1}) $. With our choice, isoperimetric theorems apply more quicquely.
Fan and Losonczy  \cite{fl} introduced this notion in $\Z^n$    in connection with an old problem of Wakeford related to
canonical forms for symmetric tensors.

 The number of distinct matching from $B$ onto $A$ will be denoted by $\mu(B,A).$

Let us define a {\em prime} group as a group having no proper finite subgroup $H .$
By elementary Group Theory, a group ia a prime group if and only if it is a torsion free group or if it has a prime order.

 The next two results are due  Losonczy in the abelian case \cite{los} and
 to Eliahou-Lecouvey in the non-abelian case \cite{el}:

 \begin{itemize}
 \item If $A$ and $B$ be finite subsets of a prime group with the same cardinality such
that  $1\notin B,$ then  $\mu(B,A)\neq 0.$
\item If  $B$ is a finite subset of a  group  such
that  $1\notin B,$ then $\mu(B,B)\neq 0.$
 \end{itemize}

The relation $1\notin B$ is obviously a necessary condition for the existence of such a pairing.
Take   a proper finite subgroup $H ,$ $h\in H\setminus \{1\}$ and $a\notin H.$ Put $B=(H\setminus \{1\})\cup \{a\}$. Suppose that  there exists a Wakeford pairing $\phi : B \rightarrow H.$ Clearly $\phi ({a})=1$ and hence $\phi ({h})\in H.$  Thus  $\phi(h)h\in H,$
 a contradiction.

 These observations made by Losonczy in \cite{los} show that    in a non-prime group, there exist finite subsets $A$ and $B$ with the same cardinality such
that  $1\notin B,$ and  $\mu(B,A)= 0.$ such a thing can not hold in prime groups, by the result mentioned above.

The results mentioned above are proved using some standard Addition Theorems presented below:

Let $A$ and $B$ be finite subsets of a  group $G$.  Kneser's Theorem \cite{knesrdensite}
states that $|AB|\ge |A|+|B|-1$ if $G$ is abelian and if $AB$ is aperiodic. The  Scherck-Kemperman Theorem states that
$|AB|\ge |A|+|B|-1$, if $A\cap B^{-1}=\{1\}.$ This result was proved first by Scherck \cite{scherck} for abelian groups
 and by Kemperman \cite{kempcompl} for arbitrary groups. As observed
  by Eliahou-Lecouvey \cite{el}, the existence of a symmetric pairing is related
 the Scherck-Kemperman Theorem. In the non-abelian case, a result due to Olson \cite{olsonjnt}
 states that $|AB|\ge |A|+|B|-|H|,$ where $H$ is a subgroup depending on $A$ and $B$.
 This last result is related to one result proved independently and shortly before
 it by the author in \cite{hejc2}. Developments of this the last result are known as
 isoperimetric results:

 A main tool in the present work is the isoperimetric approach developed by the author
\cite{hejc2,hejc3,halgebra,hast,hkemp}. Let us present briefly this method:

The subgroup generated by a set $X$ will be denoted by $\subgp{X}.$
Let $S$ be a finite subset of a group $G$ with $1\in S$.
  The {\em $kth$--connectivity}
of $S$
 is defined  as
 $$
{\kappa _k}(S)=\min  \{|(XS)\setminus X|\   :  \ \
\infty >|X|\geq k,\ X\subset \subgp{S} \ {\rm and}\ |XS| \le  |\subgp{S}|- k\},
$$
where $\min \emptyset =|\subgp{S}|-k+1$.

 We shall say that
$S$ is a {\em Cauchy subset} if   $\kappa _1(S)= |S|
-1.$
 We shall say that
$S$ is a {\em Vosper subset} if $\kappa _2(S)\ge |S|
.$
Clearly $S$ is a Cauchy subset if and only if for every
$X\subset \subgp{S}$ with $|X|\ge 1$,
 $$|XS|\ge \min \Big(| \subgp{S}|, |X|+|S|-1\Big).$$

Also, $S$ is a Vosper subset if and only if for every
$X\subset \subgp{S}$ with $|X|\ge 2$,
 $$|XS|\ge \min \Big(| \subgp{S}|-1, |X|+|S|\Big).$$

A non-empty subset $S$ of a group will be called a Chowla subset if the order of every element of $S$ is $\ge |S|+1.$
The notion  of a Chowla subset, introduced by the author in \cite{halgebra} as a relaxation of Chowla's condition in cyclic groups, allows to extend additive properties of prime groups to a large class of subsets   of an arbitrary  group.

Pairings existence is related to the Cauchy property and the Vosper's property
allows to give a lower bound for the number of distinct pairings. We shall investigate pairing
 of a Chowla subset $B$  onto an arbitrary $A$ and  obtain a lower bound for the number of distinct pairings in this case.

Recall the following  notion used by    K\'arolyi in his generalization of Vosper's Theorem \cite{karolyi}:

For a group $G,$ put $$p(G)=\min \{ |M| : M \  \mbox{is a finite subgroup of } G, \ \mbox{ with}\  2\le |M|<|G|\},$$ where
$\min \emptyset =\infty .$ In particular $p(G)=\infty$ if $G$ is a prime group.

The organization of the paper is the following:

Section 2 contains essentially known results. Section 2.1 presents K\"onig-Hall's Theorem and its
particular formulation in the Wakeford graph. Section 2.2 the Erd\H{o}s-Heilbronn averaging argument
used to give a bound for the maximum degree of the Wakeford graph ${\cal R}=\{(x,y)\in B\times A \ |\ xy \notin A\}$. Section 2.3 presents some isoperimetric formalism. In Section 3, we  prove an
inverse theorem for cofinite sets. As an application, we obtain the following result:

 If  $ B$  be a finite  Chowla subset of  a   group $G$
and  $A$ is  a  subset of  $G$ with $|A|=|B|,$ then $\mu(B,A)\neq 0.$

In particular,
 $\mu(B,A)\neq 0,$ if   $ |B|<p(G),$
where  $ B$  is a  finite subset of  a   group $G\setminus \{1\}$
 and  $A$ is a  subset of  $G$ with $|A|=|B|.$

Applied to groups with a prime order, the last result reduces to a result of
Losonczy \cite{los}. Also in the case of torsion free groups, it reduces to a result of Eliahou-Lecouvey \cite{el} generalizing
 results of Wakeford \cite{wak} and Losonczy \cite{los} in the abelian case \cite{wak}.

Assuming that $ B$ is a finite Chowla subset of  a   group $G$ and that  $A$ be a  subset of  $G$ with $|A|=|B|,$  we  show in  Section 4 that   one of the following holds:
\begin{itemize}
\item[(i)]  $\mu(B,A)\ge \min \{\frac{||B|+1}{3},\frac{|B|(q-|B|-1)}{2q-|B|-4}\},$
  \item[(ii)] For some $a\in A,$ $|Aa^{-1}\cap B|=|B|-1,$
  \item[(iii)] For some $a\in A,$ $Aa^{-1}$ is a progression,
  \end{itemize}
where $q$ denotes the cardinality of the subgroup generated by $B$.
In particular
either $\mu(B,B)
 \ge \min \{\frac{||B|+1}{3},\frac{|B|(q-|B|-1)}{2q-|B|-4}\},$ or
  for some $a\in B,$ $Ba^{-1}$ is a progression.

  Let us now construct two basic examples having a unique Wakeford pairing, where exactly one of the conditions (ii) and (iii) is satisfied:

Put $P=\{1, r, \cdots , r^j\}.$  There is clearly a unique matching from $rP$ onto $P$, where
$\phi :r^s \mapsto r^{j+1-s}.$ Thus $\mu (rP,P)=1.$

Put $P=\{1, r^2, \cdots , r^{j+1}\}$ and take $a\notin P\cup \{ r^{j+1},r^{j+2}\}.$
Set $Q=\{r^2, \cdots , r^{j+1},a\}$
There is a unique matching
$\phi :Q\rightarrow P$ with $\phi (a)=1.$  Thus $\mu (Q,P)=1.$

The last example shows that one may have $\mu(A,B)=1$ for a set  $A$ which is not a right progression.

\section{Preliminaries}

\subsection{Matchings}

Let ${\cal R}\subset V\times W$ be a relation. A {\em selection} of ${\cal R}$
is a mapping ${\cal S}:  V \rightarrow W $ such that $(x,{\cal S}(x))\in {\cal R},$ for every $x\in V.$
We shall write

$\Delta ({\cal R})=\max \{ |{\cal R}(x)|: x\in V\}$
and $\Delta ^{-1}({\cal R})=\max \{ |{\cal R}^{-1}(x)| :  x\in W\}$.

Suppose that $|V|=|W|$. A bijective selection of ${\cal R}$ is called a  {\em  matching} of ${\cal R}.$  The number of  matchings of ${\cal R}$ will be denoted by  $\mu ({\cal R}).$

We shall use the next two known results:
\begin{theorem} \label{khall} (K\"onig-Hall's Theorem) \cite{lovasz}
Let ${\cal R}\subset V\times W$ be a relation
with $|V|=|W|$. Then   ${\cal R}$ has a  matching  if and only if $|{\cal R}(Y)|\ge |Y|,$ for every  subset $Y$ of $V.$
\end{theorem}

An easy consequence of K\"onig-Hall's Theorem is the following result observed by
Hetyei, c.f. \cite{lovasz}:
\begin{corollary} \label{hetei} Let ${\cal R}\subset V\times W$ be a relation
with $|V|=|W|$
 such that $|{\cal R}(Y)|\ge |Y|+1,$  for every proper subset $Y$ of $V$.
Then   every  arc $(x,y)\in {\cal R}$ lies in some  matching of ${\cal R}$. In particular,
$\mu({\cal R})\ge \max \{\Delta ({\cal R}), \Delta ^{-1}({\cal R})\}.$
\end{corollary}

Some applications of Matching Theory may be found in the book of Lov\'asz and Plummer \cite{lovasz}.
We mention also  applications  of Theorem \ref{khall} and Corollary \ref{hetei}
by Fournier to  planar tailings. \cite{fournier}.


Let $A$ and $B$ be finite subsets of a group $G.$  The relation ${\cal R}=\{(x,y)\in B\times A \ |\ xy \notin A\}$ will be called a
{\em Wakeford graph}.
Clearly ${\cal R}(x)=A \cap (  x^{-1}\oa),$ for every $x\in B$ and ${\cal R}(X)=A\cap(X^{-1}\oa  ),$
for every $X\subset B.$

Clearly a {Wakeford pairing } from $B$ onto $A$ is just a   matching of the relation
${\cal R}.$
Thus  $$\mu (B,A)=\mu ({\cal R}).$$

We shall say that the couple $(B,A)$ is {\em matchable} if
there is a matching from $B$ onto $A$.

Notice that $xy\notin A$ if and only if
$xya\notin Aa.$ In particular
\begin{equation}\mu(B,A)=\mu(B,Aa) \label{trans}\end{equation}

For every $x\in B,$ we have ${\cal R}(x)=A\setminus (x^{-1}A).$
In particular,
\begin{equation}|{\cal R}(x)|=|A\setminus (x^{-1}A)|=|xA\setminus A|=\lambda _A(x). \label{deg}\end{equation}


Since  ${\cal R}(X)=(X^{-1}\oa)\setminus \oa =((\tilde{X}^{-1})\oa)\setminus \oa,$
 Hall's condition in this graph takes the following form:

For every $X\subset B,$ $  |X|\le |(\tilde{X})^{-1}\oa)\setminus \oa|.$
In particular, we have by K\"onig-Hall's Theorem and Corollary \ref{hetei}:

\begin{lemma} \label{KHC}
Let  $ B$  be a finite subset of  a  group $G$
 with
$1\notin B$ and let $A$ be a  subset of  $G$ with $|A|=|B|.$ Let ${\cal R}\subset B\times A$ be the Wakeford graph.
 Then  \begin{itemize}
\item $\mu (B,A)\neq 0$ if and only if for every $X\subset B,$ $  |X|\le |(\oaa\tilde{X})\setminus \oaa|.$
\item If $\mu (B,A)<\Delta ({\cal R}),$ then there exists a proper subset $X\subset B,$ with $  |X|\ge |(\oaa\tilde{X})\setminus \oaa|.$
       \end{itemize}
 \end{lemma}

\subsection{The Erd\H{o}s-Heilbronn averaging argument}

Let $G$ be a group, $B\subset G$ and  $x\in G$. The
Erd\H{o}s-Heilbronn function $\lambda$ is defined by the relation
$$
\lambda_B(x) = |(Bx)\setminus {B}|.
$$
Erd\H{o}s and  Heilbronn introduced this function in the abelian case and proved its sub-modularity \cite{EH}.
Olson generalized it to the non-abelian case. We need two  properties of this function:

\begin{lemma}[Olson \cite{olsonaa}] Let $B$ and $C$ be nonempty subsets of
a group $G$ such that $1\not\in C$. Then,
  \begin{eqnarray}
  \lambda_B(x)+\lambda_B(y)&\geq&\lambda_B(xy) .\label{eq:x+y}\\
 \sum _{x\in C} \lambda_B(x)&\geq& |B|(|C|-|B|+1). \label{eq:clique}
  \end{eqnarray}
\end{lemma}

The following lemma will be used later to show that $\Delta ({\cal R})$ is not small.

\begin{lemma} \label{EHO}
Let  $ S$ and $T$  be  finite subsets of  a   group $G$
 with
$1\notin S$ and put $q=|\subgp{S}|.$
 Then  there is an $x\in S$ such that $$\lambda _T(x) \ge \min \{\frac{|T|(|S|+\kappa _2-|T|+1)}{|S|+2\kappa _2},\frac{|T|(q-|T|-1)}{2q-|S|-4}\},$$ where $\kappa _2=\kappa _2(\tilde{S}).$

 In particular, if $\kappa _2=|S|=|T|,$
  then there is an $x\in S$ such that \begin{equation}\lambda _T(x) \ge \min \{\frac{|S|+1}{3},\frac{|S|(q-|S|-1)}{2q-|S|-4}\}.
  \label{ehol}\end{equation}
\end{lemma}

\begin{proof}
Put $\alpha =\max \{\lambda _T (x): x\in S\}.$

By the definition of $\kappa _2$, we have $|\tilde{S}^2|\ge \min (q-1,1+|S|+\kappa _2).$
Take a subset $C$ of $\tilde{S}^2$ such that
$1\notin C,$ $S\subset C$ and
 $|C|=\min (q-2,|S|+\kappa _2)$

By ( \ref{eq:clique}) and ( \ref{eq:x+y}),
\begin{eqnarray*} |T|(|C|-|T|+1)&\le& \sum_{x\in C}\lambda_T(x)\\
&=&  \sum_{x\in S}\lambda_T(x) +\sum_{x\in C\setminus S}\lambda_T(x)\\
&\le& \alpha |S|+2\alpha (|C|-|S|)\\
&=& \alpha (2|C|-|S|).\end{eqnarray*}
Thus $\alpha \ge \frac{|T|(|C|-|T|+1)}{2|C|-|S|}.$

Assume first that $|S|+\kappa _2\le q-2.$ Then $|C|=|S|+\kappa _2$ and hence
$\alpha \ge \frac{|T|(|S|+\kappa _2-|T|+1)}{|S|+2\kappa _2}.$

Assume now that $|S|+\kappa _2> q-2.$ Then $|C|=q-2$ and hence
$\alpha \ge \frac{|T|(q-|T|-1)}{2q-|S|-4}.$\end{proof}

 \subsection{Isoperimetric Preliminaries}

For a subset $X\subset G,$ we shall write $\overline{X}=G\setminus X$ and
$\tilde{X}=X\cup \{1\}.$

Let $T$ and $S$ be subsets of a group $G$ with $1\in S$. We put
$T^S=G\setminus (TS)$  and $\partial _S(T)=(TS)\setminus T.$
Clearly $G=T\cup T^S\cup \partial _S(T)$ is a partition.
We shall write $\partial^{-}_S(T)=(TS^{-1})\setminus T.$
Clearly $\partial^{-}_S(T^S)\cap T=\emptyset,$ otherwise there exist $z\in T^S$ and $y\in S$ such that
$zy^{-1}=x\in T$, and hence $z=xy\in T\cup \partial _S(T),$ a contradiction. Hence
\begin{equation}\label{id0}
\partial^{-}_S(T^S)\subset \partial_S(T).\end{equation}
The last observation, used extensively in the isoperimetric method, contains a useful duality.

Recall the following result:

\begin{theorem} \label{Cay}\cite{hejc2,hast}
Let  $ S$  be a finite subset of a group  with $1\in S$. Then there is a
finite subgroup $L\neq \subgp{S}$ generated by a subset of $S$ such that
   $\kappa_1(S)=\min (|LS|-|L|,|SL|-|L|)$.

\end{theorem}

\begin{corollary} \label{cchowla}[Proposition 2.8,\cite{halgebra}]
$\kappa _1(\tilde{S})= |S|,$ for any    finite  Chowla subset
$ S.$

\end{corollary}

\begin{proof}
By Theorem \ref{Cay},  for some finite subgroup $L\neq \subgp{S},$
 generated by a nonempty subset of $S,$ we have $\kappa_1(\tilde{S})=\min (|L\tilde{S}|-|L|,|\tilde{S}L|-|L|)\ge |S|,$ if $|L|=1$.
Assume that $|L|\ge 2$.
Since $1\in S,$ we have $\min (|LS|,|S L|)\ge 2|L|.$
 Thus
   $\kappa_1(\tilde{S})=\min (|LS|-|L|,|SL|-|L|)\ge |L|\ge |S|+1,$  contradicting the obvious inequality $\kappa_1(\tilde{S})\le |S|.$
\end{proof}

We need also the following more precise result:
\begin{theorem} \label{vchowla}[Theorem 3.2, \cite{halgebra}]
 Let  $ S$  be a finite Chowla subset. Then  $\tilde{S}$ is either a Vosper's subset or a
 progression.
\end{theorem}
Let us mention that the notions of left and right progressions
coincide for a subset containing $1.$ For this reason we shall
formulate some results using  translate copies of sets.

\section{An inverse theorem for cofinite subsets}

Inverse Theory (including the isoperimetric approach)
deals only with finite sets. In this section, we  derive an
inverse theorem for cofinite sets.

\begin{proposition} \label{cf}
Let  $ S$  be a finite Cauchy subset of  a  group $G$
 with
$1\in S.$  Let $T$ be a cofinite subset of  $G$
with $|\ot|\le |\subgp{S}|-1 .$ Then
  $|\G_S(T)|\ge |S|-1$ and \begin{equation}
T^Sa^{-1}\subset \subgp{S},\ \mbox{ for every}\  a\in T^S. \label{cf1}
\end{equation}

  Assuming moreover that $S$ is a Vosper's subset and that   $|\G_S(T)|=|S|-1$. Then
      \begin{eqnarray}
 (T^S)S^{-1}&=& T, \mbox{and}\label{cf2} \\
  |(T^S)S^{-1}|&=&
  |T^S|+|{S}|-1.\label{cf3}
  \end{eqnarray}

\end{proposition}
\begin{proof}
Put and $Z={T}^S$ and   $H=\subgp{S}.$
Take a left-decomposition $Z=Z_1\cup \cdots \cup Z_j$ modulo $H$, where  $Z_i=z_iH\cap Z,$ for some $z_i\in Z.$
By (\ref{id0}),
 $$Z_iS^{-1}\subset Z_i\cup
\partial^{-}_S (Z_i) \subset Z_i\cup \partial _S (Z_i) \subset \ot,$$ and hence $|Z_iS^{-1}|\le |\ot|\le |H|-1.$
 By the definition of $\kappa _1,$ we have $|Z_iS^{-1}|=|(z_i^{-1})Z_iS^{-1}|\ge \min (|H|, | Z_i|+|S|-1)=| Z_i|+|S|-1.$
Therefore $|ZS^{-1}|=\sum _{1\le i \le j}|Z_iS^{-1}|\ge | Z|+j(|S|-1).$
Thus we have by (\ref{id0}), $ j(|S|-1)\le |\partial^{-}_S(Z)|\le |\partial_S(T)|\le |S|-1.$
Thus $j=1,$  proving (\ref{cf1}). We have also
\begin{equation}\label{idff}
|S|-1\le |\partial^{-}_S(T^S)|\le |\partial_S(T)|
\end{equation}

Assume now  that $S$ is a Vosper's subset and that   $|\G_S(T)|\le |S|-1$.
By (\ref{idff}), $|S|-1\le |\partial^{-}_S(T^S)|\le |\partial_S(T)|\le |S|-1.$

In particular, $|(T^S)S^{-1}|=|T^S|+|\partial^{-}_S(T^S)|=|(T^S)|+|S|-1,$
proving (\ref{cf3}).

By (\ref{idff}) and (\ref{id0}),
$$
(T^S)S^{-1}=T^S \cup \partial^{-}_S(T^S)
= T^S\cup \partial _S(T)= \ot.
$$
\end{proof}

We can now show that Chowla subsets behave nicely with respect to matchability.
\begin{theorem} \label{K1}
Let  $ B$  be a finite  Chowla subset of  a   group $G.$
If $A$ be a  subset of  $G$ with $|A|=|B|,$ then $\mu(B,A)\neq 0.$
\end{theorem}

\begin{proof}
Let $X$ be a an arbitrary subset of $B.$ Put $H=\subgp{X}$ and $\tilde{X}=U.$
Clearly $U$ is a Chowla subset. By Corollary \ref{cchowla}, $U$ is a Cauchy subset.
Put $V=(\overline{A})^{-1}$ and $W=V^U.$

By Proposition \ref{cf}, $|\oaa\tilde{X}\setminus \oaa|=|\G_U (V)|\ge |U|-1=|X|.$ By Lemma \ref{KHC}, $\mu(B,A)\neq 0.$
\end{proof}

In particular,
\begin{corollary} \label{karolyit}

Let  $ B$  be a  finite subset of  a   group $G$
 with
$1\notin B$ and let $A$ be a  subset of  $G$ with $|A|=|B|.$

If  $ |B|<p(G),$ then
 $\mu(B,A)\neq 0.$
\end{corollary}

Applied to groups with a prime order, the last result reduces to a result of
Losonczy \cite{los}. Also in the case of torsion free groups, it reduces to a result of Eliahou-Lecouvey \cite{el} generalizing to the non-abelian case
 results of Wakeford \cite{wak} and Losonczy \cite{los}.

\section{ Distinct matchings}

The Vosper's property implies a bound the number of distinct matchings. We
shall illustrate this relation in the easier case of Chowla subsets.

\begin{theorem} \label{mcp}
Let  $ B$  be a finite Chowla subset of  a   group $G$ let $A$ be a  subset of  $G$ with $|A|=|B|.$
Let ${\cal R}\subset B\times A$ be the Wakeford graph.
Then  one of the following holds:
\begin{itemize}
\item[(i)]  $\mu(B,A)\ge \max \{\Delta ({\cal R}), \Delta ^{-1}({\cal R})\}.$
  \item[(ii)] For some $a\in A,$ $|Aa^{-1}\cap B|=|B|-1.$
  \item[(iii)] For some $a\in A,$ $Aa^{-1}$ is a progression.
  \end{itemize}
\end{theorem}

\begin{proof}
Suppose that (i) is not satisfied. By Lemma \ref{KHC}, for some proper subset $X$ of $B,$
we have $|(\oaa\tilde{X})\setminus \oaa|=|\G_U (V)|\le |U|-1=|X|,$ where   $\tilde{X}=U$
and $V=\oaa.$

Put $H=\subgp{X}$ and $W=V^U.$ Take $w\in W$.

By Proposition \ref{cf}, $|\G_U (V)|= |U|-1=|X|.$  By (\ref{cf2}),

$$ A^{-1}=WU^{-1}.$$

 Assume first that $|W|=1,$ and hence $W=\{w\}$. Notice that $w\in A^{-1}.$ We have clearly
 $Aw=U$. Thus $|A|=|U|=|X|+1$ and hence $|X|=|B|-1.$ Therefore $|Aw\cap B|=|X|=|B|-1,$ and (ii) holds. Assume
 $|W|\ge 2.$ By Theorem \ref{vchowla}, we have one of the two possibilities:
 \begin{itemize}
   \item $U$ is an $r$--progression, for some $r$.

   By (\ref{cf3}), $|UW^{-1}|=|W|+|U|-1=|A|\le |H|-1.$
   It follows that $W^{-1}w$ is a progression and hence $UW^{-1}=Aw$ is a progression. Thus (iii) holds.

   \item $U$ is a Vosper subset.

   Then clearly $|X|\ge |U|-1\ge 2,$ (since a subset of size $2$is not a Vosper's subset). Choose $r\in X.$

   By (\ref{cf3}), $|UW^{-1}|=|W|+|U|-1=|A|\le |H|-1.$ Since $U$ is a Vosper's subset of $H$,
    we must have $|A|=|UW^{-1}|= |H|-1.$ By (\ref{cf1}),  $A \subset Hw,$ and hence
    $|(Aw^{-1})\cap H|= |H|-1.$ But $|\subgp{r}|\ge |A|+1=|H|,$ and hence $H=\subgp{r}.$
    Since $H$ is cyclic, any subset of $H$ with size $|H|-1$ is a progression. Thus $Aw^{-1}$ is a progression and (iii) holds
   \end{itemize}
 \end{proof}

\begin{corollary} 

Let  $ B$  be a  Chowla subset with $q=|\subgp{B}|.$ Then
one of the following conditions holds:

\begin{itemize}
\item[(i)] $\mu(B,B)
 \ge \min \{\frac{||B|+1}{3},\frac{|B|(q-|B|-1)}{2q-|B|-4}\}.$
\item[(ii)]  For some $a\in B,$ $Ba^{-1}$ is a progression.
  \end{itemize}
\end{corollary}
 \begin{proof}
Let ${\cal R}\subset B\times A$ be the Wakeford graph.
Suppose that (i) and (ii) are  false. By (\ref{ehol}) and (\ref{deg}), $\Delta ({\cal R})\ge \min \{\frac{|B|+1}{3},\frac{|B|(q-|B|-1)}{2q-|B|-4}\}> \mu ({\cal R}).$  By
Theorem \ref{mcp}, for some $a\in B,$ $|Ba\cap B|=|B|-1.$

Since $a\neq 1,$ we have $|B\{1,a\}|=|B|+1.$ Since $|B|$ is less than the order of $a,$ $B$ is an $a$--progression. Thus $Ba^{-1}$ is a progression, a contradiction.
 \end{proof}


\end{document}